\documentclass[a4paper, 12pt]{scrartcl}

\usepackage{amsthm}
\usepackage{tikz-cd}
\usepackage[inline]{enumitem}
\usepackage{mathtools} 
\usepackage[utf8]{inputenc}
\usepackage[french,english]{babel}
\usepackage[T1]{fontenc}
\usepackage{epigraph}
\usepackage{hyperref}

\begin{document}

\title{Univalent Foundations and the UniMath Library}
\subtitle{The Architecture of Mathematics}
\date{}
\maketitle

\parskip = 2mm
\begin{center}

\vspace{3mm}

{\large\bfseries Anthony Bordg}\footnote{University of Cambridge, Department of Computer Science and Technology, 15 JJ Thomson Avenue, Cambridge CB3 0FD, UK. website$\colon$ \url{https://sites.google.com/site/anthonybordg/home}}$\,$\footnote{Work on this paper was supported by grant GA CR P201/12/G028}
\vspace {3mm}

{}  
\end{center}

We give a concise presentation of the Univalent Foundations of mathematics outlining the main ideas, followed by a discussion of the UniMath library of formalized mathematics implementing the ideas of the Univalent Foundations (section~1), and the challenges one faces in attempting to design a large-scale library of formalized mathematics (section~2). This leads us to a general discussion about the links between architecture and mathematics where a meeting of minds is revealed between architects and mathematicians (section~3). On the way our odyssey from the foundations to the ``horizon'' of mathematics will lead us to meet the mathematicians David Hilbert and Nicolas Bourbaki as well as the architect Christopher Alexander. 

\section{The Univalent Foundations and the UniMath Library}

\subsection{The Univalent Foundations of Mathematics}
\label{UF}

The \emph{Univalent Foundations}\cite{VVPaulBernays} of mathematics designed by Vladimir Voevodsky builds upon Martin-Löf type theory\cite{MartinLof82}, a logical system for constructive mathematics with nice computational properties that makes mathematics amenable to proof-checking by computers (i.e. by a piece of software called a proof assistant). Certified or type-checked proofs should not be mistaken for automated proofs. Even if proof assistants come with various levels of automation, either built-in for elementary steps or user-defined via the so-called tactics for less basic steps, the proof assistant only checks that man-made proofs written with it are correct.

\subsubsection{The Univalence Axiom}
\label{UA}

The main characters in Martin-Löf type theory (MLTT for short) are types and elements of these types. If $T$ is a type, then the expression $t\colon T$ denotes that $t$ is an element of $T$. In particular, if $T$ is a type and $t$, $t'$ are elements of $T$ there is a new type called the identity type of $t$ and $t'$ denoted $t =_T t'$. Sometimes for convenience we will omit the type information and we will simply write $t = t'$. When one considers only a single element $t$, i.e. $t'$ is definitionaly equal to $t$, the identity type $t =_T t$ has always at least one element denoted $\mathrm{idpath}\,t$, i.e. the expression $\mathrm{idpath}\,t\colon t =_T t$ is well-formed. This term $\mathrm{idpath}$ is called a constructor and the identity types belong to a particular class of types called \emph{inductive types}. Indeed, besides its constructors (an inductive type can have either a single constructor or many constructors), a family of types defined inductively (like the identity types are when introduced formally) obeys an induction principle. In the case of identity types, this induction principle states that given a type $T$, an element $t\colon T$, a family $F$ of types indexed by an element $t_0\colon T$ and an element $p_0\colon t =_T t_0$, if there is an element $f\colon F\,t\,(\mathrm{idpath}\,t)$ (the family $F$ instantiated with the terms $t$ and $\mathrm{idpath}\,t$), then for any elements $t'\colon T$, $p\colon t =_T t' $ there is an element of the type $F\,t'\,p$, and moreover this element is $f$ itself when $t'$ and $p$ are definitionaly equal to $t$ and $\mathrm{idpath}\,t$, respectively.  Of course, one can iterate the process of building identity types, namely given $p$ and $q$ two elements of the identity type $t =_T t'$, one can form the identity type $p =_{t =_{T} t'} q$ and so on. As it happens, these identity types lead to a very rich mathematical structure and there is a surprising connection between homotopy theory and MLTT (the latter being also coined Martin-Löf \emph{dependent type theory} in reference to these dependent types, i.e. dependent on previous types for their definition which may be inductive, like in the case of identity types, or not). Roughly, one can think of $T$ as a space, two elements $t$ and $t'$ of $T$ as points of this space, two elements $p$ and $q$ of the type $t =_T t'$ as paths from $t$ to $t'$ in the space $T$, and the elements of $p =_{t =_T t'} q$ as homotopies between the paths $p$ and $q$ and so on (the elements of the successive iterated identity types being higher homotopies). Under this correspondence $\mathrm{idpath}\,t$ is the identity path between a point $t$ and itself in the given space. Each type bearing the structure of a weak $\infty$-groupoid obtained from the tower of iterated identity types over that type. Moreover, given two types $A$ and $B$ there is also a new type denoted $A\rightarrow B$ for the type of functions between $A$ and $B$. Among these functions some of them have a distinctive property, namely their homotopy fibers\footnote{The definition of the homotopy fibers of a map is given later in \ref{subsubsec:foundations}.} are contractible\footnote{The fundamental concept of contractibility is defined later in \ref{sec:hl}.}, and they are called weak equivalences\label{weq}. Again, one forms a new type for the weak equivalences between two types $A$ and $B$, denoted $A\simeq B$. Voevodsky found an interpretation of the rules of MLTT using Kan simplicial sets where an additional axiom, the so-called \emph{Univalence Axiom}, is satisfied. The Univalence Axiom (UA for short) states a property of a universe type $U$ (interpreted as the base of a universal Kan fibration), a type whose elements are themselves types called "small" types. More specifically, first note that given two small types $A$ and $B$, by applying the induction principle of identity types (take $T:= U$, $t:= A$, and the family $F$ such that $F\,B\,p_0$ is $A \simeq B$ in the statement of the induction principle above) one defines a function \emph{eqweqmap\,A\,B}, from $A =_U B$ to $A\simeq B$, that maps the identity path to the identity equivalence when $B$ is definitionally equal to $A$. The Univalence Axiom states that for any two small types $A$ and $B$ the function ($\mathrm{eqweqmap\,A\,B}$) is a weak equivalence, giving the correct notion of equality (or path under the connection alluded to above) in the universe.

\subsubsection{The Homotopy Levels}
\label{sec:hl}

Note that in the function type $A\rightarrow B$ introduced above the type $B$ does not depend on the type $A$. Now, we can replace the type $B$ by a family of (small) types indexed by the type $A$, namely an element $F$ of type $A \rightarrow U$ (where $U$ is a universe), in this case we get a new type, the cartesian product of the family of types $F$, denoted $\displaystyle \prod_{x:A}F\,x$. Given two elements $f,g\colon \displaystyle \prod_{x:A}F\,x$, we could also ask if there is an equivalence between the identity type $f = g$ and the dependent product $\displaystyle \prod_{x:A} (f(x) = g(x))$. This equivalence (or rather the non-obvious implication) is known as \emph{function extensionality}\label{funext} and it does not hold in MLTT. Fortunately, UA does imply function extensionality, \textit{i.e.} given $A\colon U$, $F\colon A\rightarrow U$ and $f,g\colon \displaystyle \prod_{x:A}F\,x$, using UA one produces a term of the type $(\displaystyle \prod_{x:A} f(x) = g(x)) \simeq (f = g)$. Thus, the Univalence Axiom can be seen as a strong form of extensionality and the Univalent Foundations are a powerful and elegant way to achieve extensional concepts in Martin-Löf dependent type theory.\\  
Without surprise another very important type is the type of natural numbers denoted $\mathrm{nat}$. This is a second example of an inductive type. The type $\mathrm{nat}$ has two constructors, $0$ of type $\mathrm{nat}$ and $s$ of type $\mathrm{nat}\rightarrow \mathrm{nat}$ that corresponds to the successor function. The induction principle of $\mathrm{nat}$ is what one expects, namely an element of the type $$\displaystyle \prod_{P:\mathrm{nat}\rightarrow U} P\,0 \rightarrow (\prod_{n:\mathrm{nat}} P\,n \rightarrow P\,(s\,n)) \rightarrow (\prod_{n:\mathrm{nat}} P\,n)\,.$$ Finally, we would like to introduce an additional dependent type called the dependent sum type. Given a type $A$ and an element $B\colon A\rightarrow U$, we form the type of dependent pairs $(x,y)$ with $x\colon A$ and $y\colon B\,x$, denoted $\displaystyle \sum_{x:A} B\,x$. Given a small type $A$, the type $A$ might have the property that it has an element $\mathrm{cntr}\colon A$ together with for every element $x\colon A$ a path from $x$ to $\mathrm{cntr}$, \textit{i.e.} an element of $\displaystyle \prod_{x:A} x =_A \mathrm{cntr}$. The dependent sum allows us to form the type of such elements, namely $\displaystyle \sum_{\mathrm{cntr}:A} \prod_{x:A} (x =_A\mathrm{cntr})$, shortened to {\em iscontr\,A}, that corresponds to the type of proofs that $A$ seen as a space is contractible. If this last type is inhabited, \textit{i.e.} if it has an element, the type $A$ is said to be {\em contractible} and $\mathrm{cntr}$ is called {\em a center of contraction}. We are now equipped with all the tools we need to introduce the very important concept of \emph{homotopy levels}, the so-called \emph{h-levels}, that intuitively capture the fact that at some point in the tower connected with a type the iterated identity types might be contractible. First, we need to know that one is allowed to define functions over inductive types, in particular over the type of natural numbers $\mathrm{nat}$. Hence, we will define an element denoted {\em isofhlevel} of type $\mathrm{nat} \rightarrow U \rightarrow U$. To achieve this, it is enough to define $\mathrm{isofhlevel}\,0\,X$ to be $\mathrm{iscontr}\,X$ and $\mathrm{isofhlevel} \,(s\,n)\, X$ to be $\displaystyle \prod_{x:X} \prod_{y:X} \mathrm{isofhlevel} \,n \,(x =_X y)$, where $X$ is a small type. Given a small type $X$ and a natural number $n$, if the type $\mathrm{isofhlevel}\,n\,X$ is inhabited, then one says  that $X$ is of h-level $n$. The type of all types of h-level $n$ is $\displaystyle \sum_{X:U}\mathrm{isofhlevel}\,n\,X$ \footnote{This type is small with respect to a higher universe. This technical detail is unimportant for people unfamiliar with type theory.}. The types of h-level 1 are called \emph{propositions}, they are the types in which any two elements are equal. The types of h-level 2 are called \emph{sets}. For $n\geq 3$ the types of h-level $n$ are higher analogs of sets\footnote{This analogy is explained in the quote from Voevodsky that ends the current section \ref{VVquote}.}. It is possible to prove for instance that given a type $X$ and an element $n\colon \mathrm{nat}$ the type $\mathrm{isofhlevel}\,n\,X$ is a proposition, that the type $\mathrm{nat}$ is a set, or that the type $\displaystyle \sum_{X:U}\mathrm{isofhlevel}\,n\,X$ is of h-level $n+1$. Moreover, the Univalence Axiom is consistent with respect to the Law of Excluded Middle for propositions and the Axiom of Choice for sets, hence not diminishing our ability to reason about propositions or sets but increasing our ability to work with higher analogs of sets. \\
Informed by homotopy theory, the main merits of the Univalent Foundations are the realization that types in MLTT are interpreted by homotopy types (topological spaces up to weak homotopy equivalences), their corresponding stratification according to the h-levels, and the ability that types give us to build weak higher groupoids through the tower of their iterated identity types. Moreover, the Univalence Axiom gives us the ability to reason formally about structures on these higher groupoids by enforcing an equivalence principle that makes two equivalent types indistinguishable in the Univalent Foundations. Indeed, let $U_0, U_1$ be two universes with $U_0\colon U_1$ and $U_0$ being univalent. Given any family $P\colon X \rightarrow U_1$, there exists two terms $\mathrm{transportf^{P}}\colon (x =_X y) \rightarrow P x \rightarrow P y$ and $\mathrm{transportb^{P}}\colon (x =_X y) \rightarrow P y \rightarrow P x$. In particular, if one takes $U_0$ for $X$, then using the univalence axiom for $U_0$ one derives two terms of types $(A \simeq B) \rightarrow P A \rightarrow P B$ and $(A \simeq B) \rightarrow P B \rightarrow P A$, respectively. \\
The Univalent Foundations realizes the following vision of Voevodsky$\colon$
\begin{quote}
\label{VVquote}
First note that we can stratify mathematical constructions by their “level”.
There is element-level mathematics - the study of element-level objects such as numbers,
polynomials or various series. Then one has set level mathematics - the study of sets with
structures such as groups, rings etc. which are invariant under isomorphisms. The next level
is traditionally called category-level, but this is misleading. A collection of set-level objects
naturally forms a groupoid since only isomorphisms are intrinsic to the objects one considers,
while more general morphisms can often be defined in a variety of ways. Thus the next level
after the set-level is the groupoid-level - the study of properties of groupoids with structures
which are invariant under the equivalences of groupoids. From this perspective a category
is an example of a groupoid with structure which is rather similar to a partial ordering on a
set.
Extending this stratification we may further consider 2-groupoids with structures, n-groupoids
with structures and $\infty$-groupoids with structures. Thus a proper language for formalization
of mathematics should allow one to directly build and study groupoids of various levels and
structures on them.
A major advantage of this point of view is that unlike $\infty$-categories, which can be defined in
many substantially different ways the world of $\infty$-groupoids is determined by Grothendieck
correspondence, which asserts that $\infty$-groupoids are “the same” as homotopy
types. Combining this correspondence with the previous considerations we come to the view
that not only homotopy theory but the whole of mathematics is the study of structures on
homotopy types.\cite{VVExperimentalLibrary}.
\end{quote}

\subsection{The UniMath Library}

Nowadays a fraction of the community working on the Univalent Foundations is involved in the design of a library of mathematics based on the Univalent Foundations, the UniMath\cite{UniMath} project, using the Coq proof assistant. \\
The UniMath library was born in 2014 by merging three previous repositories, the repository {\em Foundations} written by Vladimir Voevodsky, the repository {\em rezk\_completion} of Ahrens, Kapulkin and Shulman, and the repository {\em Ktheory} by Daniel Grayson. As of 6 September 2018, the library has the following packages (arranged in alphabetical order)$\colon$ Algebra, CategoryTheory, Combinatorics, Folds, Foundations, HomologicalAlgebra, Induction, Ktheory, MoreFoundations, NumberSystems, PAdics, RealNumbers, SubstitutionSystems, Tactics, Topology. \\
I will briefly focus on the three packages that implement the main ideas of the Univalent Foundations as outlined in the first section\ref{UF}, namely the Foundations package, the MoreFoundations package, and the CategoryTheory package. The latter grew out of the rezk\_completion repository, hence our choice of packages covers two of the three original repositories. However, we underline that we do not do justice to the contents of the library. \\

\subsubsection{The Foundations Package}
\label{subsubsec:foundations}

The Foundations package is the very core of the UniMath library. Hence, it is intended to become a very stable part of the library. For that reason this package is more or less locked. The code in that package is highly organized, and Voevodsky thought of writing it as a mathematical "textbook" about the Univalent Foundations. \\
The file {\em PartA} therein contains, among other things, the correct formulation of the type of proofs that a given function $f\colon X\rightarrow Y$ is a weak equivalence\ref{weq}, denoted {\em isweq\,f}, which is the type $\displaystyle\prod_{y:Y}\mathrm{iscontr}(\mathrm{hfiber}\,f\,y)$ \footnote{Given a small type $A$, see section \ref{sec:hl} for the definition of the type $\mathrm{iscontr\,A}$.}, where {\em hfiber\,f\,y} is the type $\displaystyle\sum_{x:X}fx = y$ of pairs consisting of an element $x\colon X$ together with a path in $Y$ from $f(x)$ to $y$. The point being that the type $\text{isweq}\, f$ is a proposition\footnote{In the UniMath library any type whose name starts by the prefix "is" is a proposition.}, \textit{i.e.} is of h-level 1. So, that type is proof-irrelevant, hence it becomes contractible, \textit{i.e.} of h-level 0, when it is inhabited. Here, proof-irrelevance means that given a proposition $P$, there exists an equality $p\colon p_1 =_P p_2$ for all elements $p_1$, $p_2$ of $P$, \textit{i.e.} for all proofs $p_1$, $p_2$ that the proposition $P$ holds. As we have seen, the property of being a weak equivalence is instrumental in formulating the Univalence Axiom. \\
This last axiom is stated in the file {\em UnivalenceAxiom}, and it consists in postulating that the type $\text{isweq}\,(\text{eqweqmap}\,A\,B)$ is inhabited for all types $A$ and $B$ of a given universe type $U$, it means that we assume that the map ($\mathrm{eqweqmap\,A\,B}$) is a weak equivalence for all elements $A$ and $B$ of the universe type $U$. We recall\ref{UA} that ($\mathrm{eqweqmap\,A\,B}$) is the canonical map from $A =_U B$, the type of equalities between $A$ and $B$, to $A \simeq B$, the type of weak equivalences between $A$ and $B$. In this last file, one also finds various consequences of the Univalence Axiom. In that file, these consequences are themselves formulated as axioms after the corresponding implications have been proved. It is an interesting choice of design. Indeed, it allows the relevant use of the command "Print Assumptions" of the Coq proof assistant in order to track the weakest axioms needed to prove a given statement. In particular, among the consequences of the Univalence Axiom proved in that file, one should mention function extensionality\ref{funext} which reassures any mathematician that pointwise equal functions are equal, a fact that any mathematician takes for granted. \\
The file {\em UnivalenceAxiom2} is devoted to the proofs that the types corresponding to the statements of these various axioms are propositions. Hence, those types are proof-irrelevant as it should be.
Last, the file {\em PartB} of the Foundations package contains the all-important formulation of the h-levels that we have covered in \ref{sec:hl}. The stratification of types according to their h-levels discriminates well behaved propositions and sets, hence it allows elegant formalizations of mathematics at the levels of sets and categories with a proof assistant.

\subsubsection{The MoreFoundations Package}

The MoreFoundations package is less fundamental than the Foundations one, but it is still very basic. There are constructions in the Foundations package that admit many variants, those variants are not necessarily interesting in themselves, at least on a purely mathematical level, but these slight variants might come in handy in a particular context. The developers of UniMath handle this discrepancy between a core mathematical idea and its various incarnations that are used as a convenient scaffolding in the process of formalization by using the MoreFoundations package as a home for those variants. It is a nice way to avoid losing the sharp structure of the Foundations package. \\
We have already underlined in \ref{sec:hl} that the Univalent Foundations (MLTT + UA) are consistent with the Law of Excluded Middle (LEM) and the Axiom of Choice (AC) if properly stated, namely if stated at the levels of propositions and sets, respectively. The type of propositions, namely $\displaystyle \sum_{X:U}\mathrm{isofhlevel}\,1\,X$, is shortened by {\em hProp} in the UniMath library. Then, the type corresponding to LEM in UniMath\footnote{See the file {\em DecidablePropositions} of the package MoreFoundations.} is roughly $\displaystyle\prod_{P:\mathrm{hProp}}P\amalg \neg P$, where $P\amalg \neg P$ denotes the coproduct of $P$ and $\neg P$. To be more specific, it is the type $\displaystyle\prod_{P:\mathrm{hProp}}\mathrm{decidable\_prop}\,P$, where $\mathrm{decidable\_prop}\,P$ is a pair consisting of the type $P\amalg \neg P$ together with a proof that it is a proposition, allowing to prove that the type-theoretic LEM is itself a proposition. \\
The type of sets, namely $\displaystyle \sum_{X:U}\mathrm{isofhlevel}\,2\,X$, is shortened by {\em hSet} in UniMath. Then, the type corresponding to AC in UniMath\footnote{See the file {\em AxiomOfChoice} in the package MoreFoundations.} is roughly 
\[
\displaystyle\prod_{X:\mathrm{hSet}}\prod_{P:X\rightarrow U}((\prod_{x:X}\mathrm{ishinh}(P\,x))\rightarrow \mathrm{ishinh}(\prod_{x:X}P\,x))\;,
\]
where given a small type $Y$, $\mathrm{ishinh}\,Y$ is the type $\displaystyle\prod_{P:hProp}((Y\rightarrow P)\rightarrow P)$. Again, the exact type is $\displaystyle\prod_{X:\mathrm{hSet}}\mathrm{ischoicebase}\,X$, where  $\mathrm{ischoicebase}\,X$ is a pair consisting of the type $\displaystyle\prod_{P:X\rightarrow U}((\prod_{x:X}\mathrm{ishinh}(P\,x))\rightarrow \mathrm{ishinh}(\prod_{x:X}P\,x))$ together with a proof that this type is a proposition, allowing to prove that our type-theoretic AC is itself a proposition. \\
In the UniMath library the type-theoretic versions of LEM and AC introduced above are not used as axioms but only as types in order to force them to be mentioned explicitely as additional assumptions whenever they are used to prove a statement. \\
Finally, in the file entitled {\em AxiomOfChoice} there is a proof that the type-theoretic AC implies the type-theoretic LEM, a result originally due to Radu Diaconescu.

\subsubsection{The CategoryTheory Package}
 
The file {\em Categories} of the package {\em CategoryTheory} gives an analog of the univalence property for categories that mimics the pattern described for the Univalence Axiom. Indeed, given a category $\mathcal{C}$ and two elements $a$ and $b$ of $\mathrm{ob}\,\mathcal{C}$, the type of objects of $\mathcal{C}$, one defines by induction a function \emph{idtoiso}, from $a =_{\mathcal{C}} b$ to the type $\mathrm{iso} \,a \,b$ of isomorphisms from $a$ to $b$, that maps the identity path to the identity morphism when $b$ is definitionaly equal to $a$. The category $\mathcal{C}$ is univalent if for any two elements $a,b\colon \mathrm{ob}\,\mathcal{C}$ the above function is a weak equivalence. As a consequence, in the Univalent Foundations all category-theoretic constructions and proofs are invariant under isomorphism of objects of a univalent category and under equivalence of univalent categories. \\
In a subdirectory entitled {\em categories} (with a lower-case c, it should not be confused with the file of the same name), one proves that the category of sets as well as many categories of structured sets (monoids, groups, rings, modules, discrete fields, all standard algebraic structures formalized in the Algebra package) are univalent\footnote{A general theorem that "isomorphism is equality" for a large class of algebraic structures (assuming the Univalence Axiom) was proven by Thierry Coquand and Nils Anders Danielsson in 2013\cite{isoiseq}. Closely related, is the formulation of the more abstract {\em Structure Identity Principle} due to Peter Aczel, see Chapter 9 of \cite{HoTTBook}.}. \\ 
Moreover, any category is equivalent to a univalent category called its Rezk completion as established in the eponymous file {\em rezk\_completion}. Note that the univalence property itself for categories is not invariant under equivalence of categories, only under isomorphism of categories.

\section{Some Challenges  to Achieve Large-Scale Libraries of Formalized Mathematics}

\subsection{Toward Massive Collaborations in Mathematics}

Formalized mathematics leads to a high level of certification and collaboration in mathematics, in the case of the UniMath library using the Coq proof assistant for the former and the open source distributed revision control system \emph{Git}, the development platform \emph{GitHub} and a Google Group {\em Univalent Mathematics}\cite{GoogleGroupUF} for the latter. Note that certification is a prerequisite for true massive collaboration in mathematics that is otherwise hardly possible, since one would need to check by hand developments by others in order to rely on them for his own proofs and developments. The advantage of certification allows the UniMath developers to focus on the quality of the code in the library. Usually, when someone submits a contribution then volunteers make remarks and suggestions to improve the formalization, sometimes several rounds of rewriting are undertaken. However, so far the formalization of mathematics has not changed by a few orders of magnitude the traditional way of collaborating between mathematicians, which typically involves at most 2 or 3 mathematicians, as illustrated by the celebrated collaboration between Hardy and Littlewood. \\
Moreover, mathematicians have blamed the formalization of mathematics for its tediousness. Regarding an influential old formal language called Automath, N.G. de Bruijn wrote in ``A survey of the project Automath''\cite{SurveyAutomath}$\colon$
\begin{quote}
	A very important thing that can be concluded from all writing experiments is the constancy of the loss factor. The loss factor expresses what we loose in shortness when translating very meticulous "ordinary" mathematics into Automath. This factor may be quite big, something like 10 or 20, but it is constant$\colon$ it does not increase if we go further in the book. It would not be too hard to push the constant factor down by efficient abbreviations. 
\end{quote}
So, de Bruijn notes two things. First, formal proofs are longer, sometimes to an inadmissible point, as measured by the loss factor. Second, the loss factor is constant and it does not increase beyond some threshold. Regarding the first point, Freek Wiedijk	in \cite{deBruijnFactor}\footnote{See also \url{http://www.cs.ru.nl/~freek/factor/}} gives interesting data for the AutoMath, Mizar, and HOL Light systems. In several cases, what Wiedijk calls the de Bruijn Factor, roughly the ratio of a formalized text to a \TeX{} encoding of its informal counterpart, is around 4. It would be interesting to have similar data for Isabelle, a system with more automation, and for a more expressive system based on a dependent type theory like Coq, these systems might compare favourably. A de Bruijn factor equal or less than 2 might be more acceptable to a mathematician, and it is certainly a goal one should strive for. In order to succeed, in addition to more efficient support for notations as pointed out by de Bruijn, we certainly need more automation to handle the most obvious and boring parts of formal proofs. In the meantime, given enough hands and eyeballs can any substantial formalization effort be made shallow enough ? We hope that the formalization of well-known mathematics, at the undergraduate or even graduate level, is parallelizable and could benefit from a divide-and-conquer approach to build comprehensive libraries that the working mathematician could use to do research-level mathematics. In this perspective, the formalization of mathematics might be more suited to massive collaborations than a project that focuses exclusively on research-level open problems like the Polymath Project\footnote{\url{https://en.wikipedia.org/wiki/Polymath_Project}}.

\subsection{The Challenge of Scalability}

With formalized mathematics one faces the challenge of \emph{scalability}. As suggested by the second point of de Bruijn, the constancy of the loss factor, the problem is not so much about the increase of de Bruijn factor with the length of a text, but about other aspects whose scalability might be problematic. \\
We will give a simple example. With the growth of the library, its index for search becomes huge. If there is no homogeneity when possible for the names of the definitions, lemmas and theorems, then it becomes very difficult for the user to check whether some item useful for his goal has already been formalized and if so to find it in the library, for instance by guessing easily its name. \\
One could hope that the tools of machine learning could offer in the near future for instance more intelligent search support for definitions, lemmas and theorems in libraries as well as some other useful automated tools. But this perspective should not prevent developers from being very careful with the design of their library of formalized mathematics to achieve something whole.

\subsection{The Foundations$\colon$ a Never-Ending Work or an Horizon}

In the 21st century something new of the same magnitude as the Bourbaki project could wait for mathematicians.
However, Bourbaki told us the following$\colon$
\begin{quote}
If formalized mathematics were as simple as the game of chess, then once
our chosen formalized language had been described there would remain
only the task of writing out our proofs in this language, just as the author
of a chess manual writes down in his notation the games he proposes to
teach, accompanied by commentaries as necessary. But the matter is
far from being as simple as that, and no great experience is necessary to
perceive that such a project is absolutely unrealizable$\colon$ the tiniest proof
at the beginning of the Theory of Sets would already require several
hundreds of signs for its complete formalization. Hence, from Book I
of this series onwards, it is imperative to condense the formalized text by
the introduction of a fairly large number of new words (called abbreviating
symbols) and additional rules of syntax (called deductive criteria). By doing
this we obtain languages which are much more manageable than the
formalized language in its strict sense. Any mathematician will agree
that these condensed languages can be considered as merely shorthand
transcriptions of the original formalized language. But we no longer
have the certainty that the passage from one of these languages to another
can be made in a purely mechanical fashion$\colon$ for to achieve this certainty
it would be necessary to complicate the rules of syntax which govern the
use of the new rules to such a point that their usefulness became illusory;
just as in algebraic calculation and in almost all forms of notation commonly
used by mathematicians, a workable instrument is preferable to one which
is theoretically more perfect but in practice far more cumbersome. (\cite{BourbaSets}, Introduction, p. 10).
\end{quote}
As it happens, the end of the 20th century gave us such unforeseen powerful theories and proof assistants, actually not so ``complicated'' as anticipated by Bourbaki, for instance under the form of the so-called Calculus of Inductive Constructions (a dependent type theory extended with various features) as embodied in Coq, equipped with notational support to handle notations even including LaTeX and unicode characters, incorporating automatic tools like tactics, and being able to automatically generate typeset documents. Contrary to Bourbaki's expectations, packaged this way these theories have rendered the formalization of mathematics more feasible. It opens new possibilities for learning and teaching mathematics\footnote{The key step towards the widespread use of formalized mathematics could be to start teaching mathematics with the help of proof assistants, not to try very hard to gain the support of the working mathematicians. Given that present-day students are the mathematicians of tomorrow, the latter could be a consequence of the former. Moreover, a mathematical education using, at least occasionally,  proof assistants could be solidly grounded in rigorous proofs and the students could benefit from it. Unfortunately, most proof assistants might be still too difficult to use except for graduate students. Note that the Isabelle proof assistant, which offers more automation, was used recently by a group of undergraduate students in Germany, under the supervision of three advisers, to formalize partly the DPRM theorem motivated by Hilbert's Tenth Problem. See their joint paper\cite{HilbertMeetsIsabelle} presented during the FLOC 2018 conference in Oxford.}, doing mathematical research, or using mathematics in industry. Thus, time may be ripe for Bourbaki's abandoned dream. However, one should keep in mind the distinction between the formalization of mathematics using proof assistants and some ultimate foundations of mathematics. The foundations of mathematics may be a never-ending work, what Bourbaki called the ``horizon''\cite{BourbaGuedj}. Hence, the importance of a third  technical challenge (in addition to massive collaborations using proof assistants, and the scalability of libraries), the \emph{migration} of libraries, for instance from a system to a more evolved system and this is why the UniMath library uses for its development only a small subset of the Coq language. Given the numerous proof assistants and libraries of formalized mathematics on the market, migration is an important issue and old code for new proof assistants should be reused as easily as possible to become the scaffolding for new achievements. The UniMath library is intended to be a whole scalable migration-friendly library of formalized mathematics with certified proofs. With respect to large-scale formalization, one very interesting aspect of Bourbaki's project consists in noticing that its members faced large-scale architectural problems well before us, since they aimed in their series of books at a rigorous, general and self-contained, reformulation of the whole of mathematics known at the time. Hence, in this regard one can learn from Bourbaki's epoch-making project, and Armand Borel's article ``Twenty-Five Years with Nicolas Bourbaki 1949-1973''\cite{Borel25years} and Pierre Cartier's article ``The Continuing Silence of Bourbaki''\cite{CartierSilence} are informative.

\section{Architecture and Mathematics}

Alexander Grothendieck was a third-generation member of Bourbaki and when reading Grothendieck's ``R\'{e}coltes et Semailles''\cite{RSGrothendieck} one can wonder why the architectural metaphor is recurrent \footnote{I will give a few examples$\colon$ ``Je me sens faire partie, quant \`a moi, de la lign\'ee des math\'ematiciens dont la vocation spontan\'ee et la joie est de construire sans cesse des maisons nouvelles. Chemin faisant, ils ne peuvent s'emp\^echer d'inventer aussi et de fa\c{c}onner au fur et \`a mesure tous les outils, ustensiles, meubles et instruments requis, tant pour construire la maison depuis les fondations jusqu'au fa\^{\i}te, que pour pourvoir en abondance les futures cuisines et les futurs ateliers, et installer la maison pour y vivre et y \^etre \`a l'aise. Pourtant, une fois tout pos\'e jusqu'au dernier ch\^eneau et au dernier tabouret, c'est rare que l'ouvrier s'attarde longuement dans ces lieux, o\`u chaque pierre et chaque chevron porte la trace de la main qui l'a travaill\'e et pos\'e. Sa place n'est pas dans la qui\'etude des univers tout faits, si accueillants et si harmonieux soient-ils - qu'ils aient \'et\'e agenc\'es par ses propres mains, ou par ceux de ses devanciers. D'autres t\^aches d\'ej\`a l'appelant sur de nouveaux chantiers, sous la pouss\'ee imp\'erieuse de besoins qu'il est peut-\^etre le seul \`a sentir clairement, ou (plus souvent encore) en devan\c{c}ant des besoins qu’il est le seul a pressentir.'' (\cite{RSGrothendieck}, 2.5 Les h\'eritiers et le b\^atisseur); and ``Comme le lecteur l'aura sans doute devin\'e, ces "th\'eories", "construites de toutes pi\`eces", ne sont autres aussi que ces "belles maisons" dont il a \'et\'e question pr\'ec\'edemment $\colon$ celles dont nous h\'eritons de nos devanciers et celles que nous sommes amen\'es \`a b\^atir de nos propres mains, \`a l'appel et \`a l'\'ecoute des choses. Et si j'ai parl\'e tant\^{o}t de l' "inventivit\'e" (ou de l'imagination) du b\^atisseur ou du forgeron, il me faudrait ajouter que ce qui en fait l'\^ame et le nerf secret, ce n’est nullement la superbe de celui qui dit : "je veux ceci, et pas cela !" et qui se compla\^{\i}t \`a d\'ecider \`a sa guise ; tel un pi\`etre architecte qui aurait ses plans tout pr\^ets en t\^ete, avant d'avoir vu et senti un terrain, et d'en avoir sond\'e les possibilit\'es et les exigences.'' (\cite{RSGrothendieck}, 2.9); and again ``C'\'etait peut-\^etre l\`a la principale raison pour laquelle les maisons que je prenais plaisir \`a construire sont rest\'ees inhabit\'ees pendant le longues ann\'ees, sauf par l'ouvrier ma\c{c}on lui-m\^eme (qui \'etait en m\^eme temps aussi l'architecte, le charpentier etc.).'' (\cite{RSGrothendieck}, 18.2.8.3 Note 135).}. Actually, there is a meeting of minds between great architects and great mathematicians linked by an abstract approach of space with surprisingly at the same time a feeling of its organic life. This abstract approach of space is  remarkable in the great architectural theoreticians, like for instance Frank Lloyd Wright, Le Corbusier, or Christopher Alexander, in the sense that one can naively believe their prime business is the 3-dimensional space embodied in a house, a building or a city, it is certainly true but it goes beyond. We notice that architects have been facing large-scale problems for long, they have been challenging them and they have offered their thoughts. Christopher Alexander in ``The Nature of Order''\cite{TNoOI} develops what he calls \emph{wholeness} to answer these challenges. Wholeness is precisely what is lacking in most libraries of formalized mathematics despite the fact it is an important feature of mathematics. Note that \emph{wholeness} in Alexander's work is a specific concept defined by 15 properties\footnote{\begin{enumerate*} \item Levels of scale \item Strong centers \item Boundaries \item Alternating repetition \item Positive space \item Good shape \item Local symmetries \item Deep interlock and ambiguity \item Contrast \item Gradients \item Roughness \item Echoes \item The void \item Simplicity and inner calm \item Not-separateness \end{enumerate*}}. We will not discuss here each of those properties, but only a few that seem more relevant with respect to mathematics, since it will be probably hopeless to search for a precise dictionary between Alexander's properties and some corresponding features of mathematics. These properties are an interesting attempt to capture what "organic" and "life" could mean for man-made artefacts like architectural works which are Alexander's main concern. In his 1900 address to mathematicians {\em Mathematical Problems} Hilbert mentioned this organic feature of mathematics in the following perceptive insights$\colon$ 
\begin{quote}
The problems mentioned are merely samples of problems, yet they will suffice to show how rich, how manifold and how extensive the mathematical science of today is, and the question is urged upon us whether mathematics is doomed to the fate of those other sciences that have split up into separate branches, whose representatives scarcely understand one another and whose connection becomes ever more loose. I do not believe this nor wish it. Mathematical science is in my opinion an indivisible whole, an organism whose vitality is conditioned upon the connection of its parts. [\ldots] We also notice that, the farther a mathematical theory is developed, the more harmoniously and uniformly does its construction proceed, and unsuspected relations are disclosed between hitherto separate branches of the science. So it happens that, with the extension of mathematics, its organic character is not lost but only manifests itself the more clearly.\cite{Hilbert23problems}.
\end{quote}
If wholeness is a feature of mathematical science according to Hilbert, Alexander regrets its absence in most of modern, dead and dull, architectural works while this property is shining in some of the great artistic works of the past. Alexander's concerns may not be widely shared by present-day architects, but they are not without resonances among other great architects as testified by the organic architecture of Frank Lloyd Wright or Tadao Ando's obsession with making light vibrant, alive, through the use of concrete.  If Alexander does not mention mathematics, a mathematician cannot help but think that mathematical entities and mathematics as a whole (``\emph{la} mathématique'' of Bourbaki, using a singular on purpose) display to a great extent this pervasive organic character, a life of their own, a wholeness. \\
Mathematical entities are like the \emph{centers} of Alexander, the elementary components of any system that make it alive, but with the subtlety that a center cannot be isolated from other centers but needs to be understood in a mutual recursive relation with other centers (\cite{TNoOI}, p.116). We can think about prime numbers in terms of what Alexander calls \emph{strong centers} (\cite{TNoOI}, p.151), centers that focus our attention and engage us, the set of primes numbers being described by Alain Connes as the heart of mathematics \cite{ConnesInterview} (compare Connes's vivid organic metaphor in the interview with the dull mechanic metaphor of the interviewers using a coffee machine). \\
In mathematics strong centers are not only displayed as specific mathematical entities but also in proofs, proving being at the core of the activity of the working mathematician. Indeed, any good mathematical proof has its own architecture. This architecture revolves around the main ideas that provide the flesh of the proof. In a given proof there are as many strong centers as there are main ideas, fitting together thanks to {\em boundaries} which are the glue of the inner workings of the mathematical mind, the hypotheses and the conclusion being the initial boundary and the last boundary, respectively. Simple proofs have usually only one center, more elaborate proofs may have many centers, but it does not matter, centers are always what make things click. One can define the strong centers in a proof as the main ideas such that handed to any mathematician with the appropriate training he will not fail to reconstruct the proof on his own. In most proof assistants, formal proofs have no structure\footnote{Some proof assistants like Isabelle have structured proofs (in the case of Isabelle thanks to an additional layer called the Isar language), but there is still room for improvement.}. It is an important issue. This is the case for instance in the Coq proof assistant, and as a consequence in the UniMath library, where a formal proof is basically a sequence of tactics lacking the structure of its informal counterpart. Even if one would not intend to read formal proofs, for instance to get pedagogical insights, but only to get certificates of correctness from them, then one still needs to maintain on a regular basis the code in a library to take into account revisions that might have been proposed. Since changes pushed in a library can break some proofs, the task of repairing broken code (and in particular broken proofs) is made harder by the lack of structure in formal proofs that makes them barely legible. \\
Also, there might be some \emph{roughness} in the sense of Alexander in the distribution of prime numbers mentioned above, roughness being an elusive property$\colon$
\begin{quote}
Things which have real life always have a certain ease, a morphological roughness. This is not an accidental property. It is not a residue of technically inferior culture, or the result of hand-craft or inaccuracy. [\ldots] It is an essential feature of living things, and has deep structural causes. [\ldots] Roughness does not seek to superimpose an arbitrary order over a design, but instead lets the larger order be relaxed, modified according to the demands and constraints which happen locally in different parts of a design. (\cite{TNoOI}, p.210).
\end{quote} 
Moreover, it suffices to quote Hilbert again to find traces in mathematics of other properties of Alexander like \emph{local symmetries}, \emph{deep interlocking} and \emph{echoes}$\colon$ 
\begin{quote}
For with all the variety of mathematical knowledge, we are still clearly conscious of the similarity of the logical devices, the relationship of the ideas in mathematics as a whole and the numerous analogies in its different departments.\cite{Hilbert23problems}.
\end{quote} 
The libraries of formalized mathematics are annoyingly lacking echoes. Often proof assistants miss some nice built-in features. For instance, these libraries have no counterpart as convenient as the index of a book, the easy search function of a PDF file or the table of contents of a book\footnote{The good practices of writing a short table of contents at the top of a file starting a new formalization and a bibliography at the end are surprisingly not even included in the style guide (\url{https://github.com/UniMath/UniMath/blob/master/UniMath/README.md}) of UniMath as of 6 September 2018.}, let alone clickable keywords for pop-up windows to remind the reader about definitions.\footnote{Again, the Isabelle theorem prover and its bundled editor jEdit, even if not perfect, have built a competitive advantage with search support and clickable keywords.} \\
So far, formalized mathematics has focused only on making impossible to write faulty proofs, in doing so it has done nothing for making proofs easier to read. Quite the contrary, formalized mathematics is much harder to read than everyday mathematics and this can explain why it has encountered considerable resistance from mathematicians. One should not forget that mathematicians spend a lot of time reading mathematics, not only doing or writing it. Formalized mathematics has forgotten the communication function of written mathematics, and this is a problem not only with respect to mathematicians but also for students, especially if one believes that teaching mathematics with proof assistants may have pedagogical value and is a necessary milestone towards the widespread use of formalized mathematics.\footnote{See also the footnote 8 on that point.} Formalized mathematics is not responsive to the reader, this strong center of the subjective experience of mathematics, that changes across the mathematical community. In this sense, it lacks the property that Alexander coined {\em gradients}, the adaptive result in design when conditions vary. A simple solution should be to have expensible/collapsible parts in proofs, so that every reader, while reading a proof, can set for himself the level of details according to his background and ability. Hopefully, this feature would allow to hide very low-level details that make reading formal proofs cumbersome. Few libraries of formalized mathematics are really designed with the reader in mind. \\ 
Finally, could it be that the last fifteenth property of Alexander, \emph{not-separateness}, the experience of ``a living whole as being at one with the world'' (\cite{TNoOI}, p.230), is the ``unreasonable effectiveness of mathematics in the natural sciences'' emphasized by Eugene Wigner  \cite{WignerUnreasonableEffectiveness} ? \footnote{I have discovered a truly remarkable answer to this question which this footer is too small to contain.} \\
In the same way an architect try to realize the unfolding in space of a form through \emph{levels of scale} (\cite{TNoOI}, p.145), \textit{i.e.} the property that consists in the presence of centers at a wide range of scales, mathematicians unfold their axioms through mathematical entities and theorems. The levels of scale are apparent for instance in the definition/lemma/theorem/corollary structure of a mathematical book or article. Both architects and mathematicians are happy when this unfolding looks like the unfolding of an organism from the seed within. The axioms of mathematics are the seeds, the labour of mathematicians are the ground that nurtures the seeds, and the mathematical entities and theorems are the resulting landscape with its wide open horizon. Some parts of this landscape are \emph{jardins à la française}, some others are English gardens, both with their respective supporters. Most parts of this landscape secretly aspire to the inner peace of Japanese gardens, natural but neither artificial nor wild. This living unfolding, from the axioms to the theorems, is the stuff mathematical objects are made from. In the case of the Univalent Foundations of mathematics the unfolding of shapes, namely types of various h-levels, a concrete example of \emph{levels of scale} relevant in our context, what appears less directly and less smoothly in sets-based mathematics as homotopy types, is remarkable. The enlarged notion of  life, of living structures, coined \emph{wholeness} by Alexander, can help to understand in particular where the platonistic attitude of mathematicians comes from. \\
As pointed out earlier, Alexander underlines the interplay of centers with the use of \emph{boundaries}(\cite{TNoOI}, p.150) which separate a center from others and at the same time unite them. For a second example, think in mathematics about the locus where two topics or two theories meet, share some methods and that could possibly merge in the future as a result. But Alexander seems to miss the point that sometimes the life of some parts may be at the expense of others. This full dynamics was noted by Hilbert$\colon$ 
\begin{quote}
[\ldots] let me point out how thoroughly it is ingrained in mathematical science that every real advance goes hand in hand with the invention of sharper tools and simpler methods which at the same time assist in understanding earlier theories and cast aside older more complicated developments. It is therefore possible for the individual investigator, when he makes these sharper tools and simpler methods his own, to find his way more easily in the various branches of mathematics than is possible in any other science.\cite{Hilbert23problems}.
\end{quote} 
The result of these sharper tools and simpler methods, won after the struggle, that ease the orientation of mathematicians in the whole of their science is the \emph{simplicity and inner calm} put forward by Alexander, and described by him as the 
\begin{quote}
quality [\ldots] which is essential to the completion of the whole. [\ldots] The quality comes about when everything unnecessary is removed.(\cite{TNoOI}, p.226).
\end{quote}
The regular clean-up and reorganizations in mathematics mentioned by Hilbert above might be the analog of evolution in the biological world and a condition for a renewal of creativity, biological systems being a paradigm of wholeness. Of course, biological evolution as understood by modern biology is a blind process, while reorganizations are made on purpose by mathematicians and some mathematicians have the platonistic feeling to be guided by an independent architectural principle of some kind, to discover rather than to invent mathematical objects. One could think this prompted Darwin to say 
\begin{quote}
I have deeply regretted that I did not proceed far enough at least to understand something of the great leading principles of mathematics, for men thus endowed seem to have an extra sense.\cite{DarwinBio}
\end{quote}
, but one could also see in this quote a reference to the complementary ability of the mathematician, like the artist, to tap into the subconscious mind as pointed out by Jacques Hadamard\cite{HadamardInvention}.

\section{Conclusion}

In this article, following Alexander's approach, we have tried to underline a few strong centers in the foundations of mathematics, namely the Univalent Foundations, the UniMath library, the Bourbaki's cathedral of mathematics, twisting some architectural threads in a mathematical landscape. However, we have only sketched the boundaries between these centers to allow for at least some wholeness in the odyssey promised in the abstract. \\
Some mathematicians are afraid that formalization could disrupt their flow of work, their inner music, and this may be indeed a real danger if one is not able to cleverly design organic libraries to allow smooth reorganizations on a regular basis. But if we are sensitive to the wholeness of our library this danger could be avoided. By facing new large-scale challenges in design formalized mathematics could offer new opportunities. The Alhambra (close to Alexander's heart) located in Granada (Spain), started in 889, still stands shadowing our mortality, in the same way can the libraries of formalized mathematics stand the test of time ?

\subsection*{Acknowledgments}

The author would like to thank Benedikt Ahrens, Thierry Coquand, and an anonymous referee for their useful comments and suggestions.

\end{document}